\documentclass[11pt]{article}
\usepackage[hidelinks]{hyperref}
\usepackage{comment}
\usepackage[margin=1in]{geometry}
\usepackage{amsmath,amssymb,amsthm,mathtools}
\usepackage{microtype}
\usepackage{todonotes}
\usepackage{authblk}
\usepackage[ruled,vlined,linesnumbered]{algorithm2e}

\newcommand{\bGam}{\bar{\Gamma}}
\newcommand{\bQ}{\bar{Q}}
\newcommand{\pluspart}[1]{\left[#1\right]^{+}}
\newcommand{\TWOSEL}{\textup{\textsf{2-Sel}}}
\newcommand{\ONESEL}{\textup{\textsf{1-Sel}}}
\newcommand{\DDID}{\textup{\textsf{DDID}}}
\newcommand{\SEL}{\textup{\textsf{Sel}}}
\newcommand{\Equal}{\textup{\textsf{EqualDev}}}

\newcommand{\PARTITION}{\textup{\textsf{Equal-Cardinality-Partition}}}

\newcommand{\sumh}{2H}
\newcommand{\Q}{\mathcal{Q}}

\newcommand{\Y}{\mathcal{Y}}

\newcommand{\obsxi}{\mu}

\newcommand{\ddid}[2]{{\textsf{DDID}-#1-#2}}

\newcommand{\R}{\mathbb R}
\newcommand{\Z}{\mathbb Z}
\newcommand{\N}{\mathbb N}

\newcommand{\set}[2]{\left\{#1 \; \left|\;\; #2 \right.\right\}}

\newcommand{\xc}[1]{{\color{black}#1}}

\newtheorem{remark}{Remark}
\newtheorem{theorem}{Theorem}
\newtheorem{lemma}{Lemma}
\newtheorem{proposition}{Proposition}
\newtheorem{corollary}{Corollary}
\newtheorem{observation}{Observation}

\newcommand{\argmax}{\operatorname*{arg\,max}}

\title{Robust decision-dependent information discovery is harder than robust optimization}

\author[1]{Xiaoyu Chen}
\author[2]{Samuel German}
\author[3]{Jérémy Omer}
\author[1]{Michael Poss}

\affil[1]{LIRMM, University of Montpellier, CNRS, Montpellier, France}
\affil[2]{Department of Computer Science and Engineering, University of California San Diego, USA}
\affil[3]{IRMAR, INSA Rennes, Rennes, France}
\date{}

\begin{document}
\maketitle

\begin{abstract}
In this paper, we consider a four-stage variant of robust optimization with cost uncertainty described by the budget uncertainty polytope, known as decision-dependent information discovery (DDID). In this problem, the decision maker first queries a subset of the cost coefficients. Then, the adversary responds by revealing worst-case values for the queried elements, after which the decision maker chooses a feasible solution while facing the remaining uncertainty on the unqueried costs. Our contribution is twofold. On the negative side, we prove that, the DDID is NP-complete even when the underlying nominal problem is the 2-selection problem, for which only a polynomial number of feasible solutions exist. \xc{This contrasts with the min-max robust counterpart with budget uncertainty, known to be not harder than the nominal problem.} On the positive side, we prove that the more general $p$-selection problem can be solved in polynomial time under the additional assumption that all deviations are equal.\\

\noindent\textbf{Keywords:} robust optimization, DDID, budget uncertainty, complexity, selection problem
\end{abstract}

\section{Introduction}

Robust optimization is a well-established framework to handle optimization problems in the presence of uncertainty, by considering that some of the parameters may take any value in a given convex set, and taking the decision that is best assuming worst-case realization of the uncertain parameters. A fundamental extension of robust optimization arose with the introduction of a preliminary level of information discovery. In this model, the decision maker can decide to query some of the uncertain parameters before having to decide on the solution to implement. The resulting problem, named robust optimization with Decision-Dependent Information Discovery, was first introduced in~\cite{VayanosGY26}, and is motivated by several applications, such as urban planning, project management, resource allocations, scheduling. In their seminal paper, the authors propose heuristic solutions based on $K$-adaptability~\cite{hanasusanto2015k}, which they assess on several numerical examples, including a case-study illustrating how preference elicitation can improve the US kidney allocation system. Despite its practical relevance, obtaining good numerical solutions to the problem, even for the special case of cost budget uncertainty, has remained very challenging. 

Introducing the query set $\Q\subseteq2^{[n]}$, the nominal feasibility set $\Y\subseteq 2^{[n]}$, the budget polytope $\Xi_\Gamma=\set{\xi\in\R^n}{\sum_{i\in[n]}\frac{\xi_i}{d_i}\leq \Gamma,  0\leq\xi\leq d}$ given by the deviations $d>0$ and budget $\Gamma> 0$, this special case can be cast as
\begin{equation}
\label{eq:DDIDopt}\tag{DDID-OPT}
\min_{Q\in \Q} \max_{\obsxi\in \Xi_\Gamma} \min_{Y\in \Y}
         \max_{\xi \in \Xi_\Gamma(Q,\obsxi)} \sum_{i\in Y}(c_i+\xi_i),
\end{equation}
where $c$ represents the nominal cost vector and $\Xi_\Gamma(Q,\obsxi)=\set{\xi\in \Xi_\Gamma}{\xi_i = \obsxi_i\text{ for each }i\in Q}$, whose defining
constraint ensures that the cost coefficients already observed cannot be altered once $Y$ has been fixed. The four stages of~\eqref{eq:DDIDopt} can be read from left to right as follows. The decision maker first selects which components of the cost vector to observe, through the measurement set $Q$; nature then reveals the values $\obsxi$ taken by the observed components; in view of this partial information, the combinatorial solution $Y$ is chosen; nature finally picks the worst realization $\xi$ that is consistent with what has been revealed.

Problem~\eqref{eq:DDIDopt} was further studied by \cite{ParadisoGDT25}, which improved over~\cite{VayanosGY26} by speeding up the $K$-adaptability approach and proposing a Benders decomposition algorithm based on combinatorial cuts. \cite{ParadisoGDT25} also proposed a detailed case study focusing on collecting medicine crates at the Alrijne Hospital, in the Netherlands. A significant step forward in the numerical solution of~\eqref{eq:DDIDopt} arose with~\cite{OmerPossRougier2024}, which provided (i) a compact mixed-integer linear programming reformulation for DDID essentially assuming that $\Y$ is easy-to-optimize, and (ii) a branch-and-price algorithm extending the latter formulation to problems for which  $\Y$ is not easy-to-optimize. These formulations made it possible to solve DDID exactly for larger problems than before.

The above papers have provided different algorithms to obtain good approximate or exact solutions to~\eqref{eq:DDIDopt}, focusing on algorithms that work well numerically. In this paper, we focus on the decision counterpart of~\eqref{eq:DDIDopt}, so given a threshold $z$, we ask whether 
\begin{equation}
\label{eq:DDID}\tag{DDID}
\min_{Q\in \Q} \max_{\obsxi\in \Xi_\Gamma} \min_{Y\in \Y}
         \max_{\xi \in \Xi_\Gamma(Q,\obsxi)} \sum_{i\in Y}(c_i+\xi_i) \leq z.
\end{equation}
The complexity of~\eqref{eq:DDID} has remained largely open since its introduction in~\cite{VayanosGY26}, contrasting with our understanding of simpler robust problems~\cite{BertsimasS03,BuchheimK18,robook}. \xc{More specifically, it is well-known~\cite{BertsimasS03} that the problem
\begin{equation}
\label{eq:minmax}
\min_{Y\in \Y}
\max_{\xi \in \Xi_\Gamma} \sum_{i\in Y}(c_i+\xi_i) \leq z
\end{equation}
is easy as soon as $\Y$ is easy-to-optimize (meaning that optimizing a linear function over $\Y$ is easy). Our purpose herein is to understand whether the positive result available for~\eqref{eq:minmax} carries over to~\eqref{eq:DDID}.} We know already from~\cite{OmerPossRougier2024} is that computing the objective function
$$
\Phi(Q)\;:=\; \max_{\obsxi\in \Xi_\Gamma} \min_{Y\in \Y}
         \max_{\xi \in \Xi_\Gamma(Q,\obsxi)} \sum_{i\in Y}(c_i+\xi_i)
$$
is easy as soon as $\Y$ is easy-to-optimize. Hence, enumerating all $Q\in\Q$ leads to a polynomial algorithm for~\eqref{eq:DDID} when $\Y$ is easy-to-optimize and $|\Q|$ is polynomially bounded. Yet, enforcing $|\Q|$ to be polynomially bounded is a very strong assumption, since natural query sets would, at the very least, involve query budgets that are part of the input~\cite{OmerPossRougier2024,ParadisoGDT25,VayanosGY26}, naturally leading to sets $\Q$ of superpolynomial cardinalities. On the other hand, if $\Q$ or $\Y$ is hard-to-optimize, it is easy to see that~\eqref{eq:DDID} is NP-hard (straightforward for $\Y$, detailed in Section~\ref{sec:reformulations} for $\Q$). Hence, the most interesting complexity problem lies in the case where both $\Q$ and $\Y$ are easy-to-optimize. The main contribution of this paper is to answer this question, showing that DDID is NP-hard even when $\Q=\set{Q\subseteq[n]}{|Q|=q}$ and $\Y=\set{Y\subseteq[n]}{|Y|=2}$; observe that our result thus involves a particularly simple set $\Y$, the cardinality of which is only $O(n^2)$! We complement this hardness result by showing that~\eqref{eq:DDID} is in NP, dismissing the possibility that higher levels of the polynomial hierarchy above the first are needed to establish an upper bound on its complexity, which possibility is the case with other adjustable variants of robust optimization problems with continuous budget uncertainty~\cite{goerigk2024complexity,poss2026complexity}. Last, we prove that two special cases are polynomially solvable. The first one involves the set $\Y=\set{Y\subseteq[n]}{|Y|=1}$ and the algorithm follows quite naturally from the preliminary observations derived for~\eqref{eq:DDID}. The second one involves set $\Y=\set{Y\subseteq[n]}{|Y|=p}$, for some integer $p>0$, and further assumes that all deviations $d$ are equal, which can arise in scheduling or measurements from identical sensors; 
the assumption is also a natural choice when little information on the individual deviations is available. The proof of this latter result relies on a more subtle case analysis. 

We mention that other works have been dedicated to variants of~\eqref{eq:DDID}. Perhaps closest to the present work are~\cite{CHEN2026277} and~\cite{chen2026decision}. On the one hand,~\cite{chen2026decision} considers the variant of~\eqref{eq:DDID} with uncertainty on a single constraint of $\Y$ (and no cost uncertainty) and proves that such problems are easy when some cardinality constrained variant of $\Y$ are easy-to-optimize and a constant number of items can be queried, NP-hard otherwise. On the other hand,~\cite{CHEN2026277} focuses on the discrete budget uncertainty set (the extreme points of $\Xi$) and $\Y$ being the selection set, and considers either cost or constraint uncertainty, providing different types of results: polynomiality of special cases, hardness for more general settings, and a MILP reformulation. Other variants of~\eqref{eq:DDID} have also been considered in the literature, such as~\cite{jin2024distributionally}, who consider a distributionally robust version, solved by combining a decomposition algorithm with $K$-adaptability, and~\cite{ley2026decision} who consider an extension called controllable uncertainty, generalizing~\eqref{eq:DDID} to the case where the observation reduces the uncertainty around a (possibly unknown) hedging point.

The rest of this paper is structured as follows. Section~\ref{sec:NP} formally defines problem~\eqref{eq:DDID} as a decision problem and proves its membership in NP. Section~\ref{sec:reformulations} introduces an adversarial reformulation of~\eqref{eq:DDID} and illustrates how the latter leads to the problem being easy for $\Q=\set{Q\subseteq[n]}{|Q|=q}$ and $\Y=\set{Y\subseteq[n]}{|Y|=1}$. Section~\ref{sec:hardness} leverages the adversarial reformulation to prove the hardness of~\eqref{eq:DDID}, using a reduction from \PARTITION. Section~\ref{sec:equaldev} turns to the problem with equal deviations, providing a polynomial time algorithm for the problem.

\section{NP membership}
\label{sec:NP}

The precise definition of~\eqref{eq:DDID} as a decision problem requires the definition of the instances and how these are related to sets $\Q$ and $\Y$. We consider in this section the general version of the problem, named $\DDID$, by assuming that both $\Q$ and $\Y$ are defined by linear constraints. Specifically, an instance $I$ of problem $\DDID$ is given by rational matrices $A$ and $B$, rational vectors $a,b,c,d$, $d>0$, and rationals $\Gamma>0$ and $z$. Given these, one defines $\Q=\set{Q\subseteq[n]}{\sum_{j\in Q}A_{\bullet j}\leq a}$ and $\Y=\set{Y\subseteq[n]}{\sum_{j\in Y}B_{\bullet j}\leq b}$, where $M_{\bullet j}$ denotes the $j$-th column of matrix $M$, and $\DDID$ asks whether there exists a solution to~\eqref{eq:DDID}, or said differently, whether there exists $Q\in\Q$ such that $\Phi(Q)\leq z$.

The membership of $\DDID$ in NP follows from folklore results from linear programming and the non-compact reformulation based on the convexification of $\Y$ proposed by~\cite[Section 5.1]{OmerPossRougier2024}. Let us first state a lemma on the existence of a small solution to any mixed-integer linear program with polynomially many integer variables and constraints, but possibly large number of continuous variables. 

\begin{lemma}
\label{eq:MILPinNP}
Let $P$ be a mixed-integer linear program with $m$ constraints and $p$ integer variables, where $m$ and $p$ are polynomially bounded in the input size and the integer variables have polynomially bounded encoding length, and assume its optimal solution is finite. Then, the optimal solution cost $z$ has polynomially bounded encoding length and $P$ admits an optimal solution with at most $m+1$ non-zero continuous components, with polynomially bounded encoding length.
\end{lemma}

\begin{proof}
Let $\hat x$ be an optimal solution to $P$ of cost $z$. Define $P'$ from $P$ by fixing all integer variables to their value in $\hat x$. Problem $P'$ is an LP of cost at most $z$, so there exists a basic feasible solution which has at most $m+1$ non-zero components and satisfies $M\,x_M=h_M$ for a non-singular submatrix $M$ of order at most $m+1$. By Cramer's rule the encoding length of each component of $x_M$ is polynomially bounded in the input size. This also proves that the encoding length of the optimal cost $z$ is polynomially bounded.
\end{proof}

\begin{corollary}
\label{cor:NP}
    $\DDID$ is in NP.
\end{corollary}
\begin{proof}
Consider the formulation (DDID-CG) proposed in~\cite[Section 5.1]{OmerPossRougier2024} for~\eqref{eq:DDIDopt}. It satisfies the two assumptions of Lemma~\ref{eq:MILPinNP}, so that any YES-instance possesses a small certificate and the problem thus belongs to NP.
\end{proof}

\section{Reformulation of \DDID}
\label{sec:reformulations}

We provide next a reformulation of~\eqref{eq:DDID} that will be instrumental in proving the hardness of the problem, as well as proving one easy special case. Introducing an additional continuous variable $z$, we obtain the epigraphical reformulation of $\Phi(Q)$ as
\begin{subequations}
\begin{align}
\Phi(Q)\;=\; \max_{\mu,z}\; &z \\
\text{s.t.}\; &z\leq \min_{Y\in\Y}\max_{\xi\in\Xi_\Gamma(Q,\mu)}\sum_{i\in Y}(c_i+\xi_i),\label{eq:genPhiminmax}\\
& \mu\in\Xi_\Gamma.
\end{align}
\end{subequations}
Introduce the complement of $Q$ as $\bar Q:=[n]\setminus Q$ (which is empty if there is $Q\in\Q$ equal to $[n]$) and $\bGam_Q(\mu) := \Gamma - \sum_{e\in Q}\frac{\mu_e}{d_e}$, and observe that 
\begin{equation}
\label{eq:splitting}
\max_{\xi\in\Xi_\Gamma(Q,\mu)}\sum_{i\in Y}(c_i+\xi_i)=\sum_{i\in Y}c_i+\sum_{i\in Y\cap Q}\mu_i+\max_{\xi\in\Xi_{\bGam_Q(\mu)}}\sum_{i\in Y\cap\bar Q}\xi_i.
\end{equation}
Linearizing the minimum over $Y$ of constraints~\eqref{eq:genPhiminmax} and using~\eqref{eq:splitting} leads to
\begin{subequations}
  \label{form:Phi-lin}
\begin{align}
\Phi(Q)\;=\; \max_{\mu,z}\; &z \\
\text{s.t.}\; &z\leq \sum_{i\in Y}c_i + \sum_{i\in Y\cap Q}\mu_i+\max_{\xi\in\Xi_{\bGam_Q(\mu)}}\sum_{i\in Y\cap\bar Q}\xi_i, && \forall Y \in\Y,\\
& \mu\in\Xi_\Gamma
\end{align}
\end{subequations}
We propose next to compute $\Phi(Q)$ through a binary search on the value of $z$. Specifically, we define $\gamma_Q(z)$ as the minimum budget value $\Gamma$ needed by the adversary so the cost of every $Y\in\Y$ is at least $z$, formally
\begin{subequations}
  \label{form:gamma}
\begin{align}
\gamma_Q(z)\;:=\; \min_{\mu,\bar \Gamma}\; &\bar\Gamma+ \sum_{i\in Q}\frac{\mu_i}{d_i} \\
\text{s.t.}\; &z\leq \sum_{i\in Y}c_i + \sum_{i\in Y\cap Q}\mu_i+\max_{\xi\in\Xi_{\bar \Gamma}}\sum_{i\in Y\cap\bar Q}\xi_i, && \forall Y \in\Y,\label{form:gamma:Y}\\
& 0\leq \mu \leq d,\\
& 0\leq \bar\Gamma \leq \bar\Gamma^{\max},
\end{align}
\end{subequations}
where we observe that the upper bound $\bar\Gamma^{\max}:=\max_{Y\in\Y}|Y|$ is satisfied by any optimal solution.
\begin{proposition}
\label{prop:Phiandgamma}
For any $z\in\R$, $\Phi(Q)< z\iff\gamma_Q(z)> \Gamma$.
\end{proposition}
\begin{proof}
Let us prove the equivalent statement $\Phi(Q)\geq z\iff\gamma_Q(z)\leq \Gamma$. Suppose $\gamma_Q(z)\leq \Gamma$ and let $(\mu^*,\bar\Gamma^*)$ form an optimal solution of~\eqref{form:gamma}, which we assume to satisfy $\mu_i=0$ for each $i\in\bar Q$ without loss of generality. Then, $\bar\Gamma^*+ \sum_{i\in Q}\frac{\mu_i^*}{d_i}\leq \Gamma$ implies that $ \bar\Gamma^*\leq \bar\Gamma_Q(\mu)$ so $\max_{\xi\in\Xi_{\bar \Gamma^*}}\sum_{i\in Y\cap\bar Q}\xi_i\leq \max_{\xi\in\Xi_{\bGam(Q,\mu)}}\sum_{i\in Y\cap\bar Q}\xi_i$, and $(\mu^*,z)$ is feasible for~\eqref{form:Phi-lin},  so $\Phi(Q)\geq z$. Conversely, assume $\Phi(Q)\geq z$ and let $(\mu^*,z^*)$ be an optimal solution to~\eqref{form:Phi-lin}. Then $(\mu^*,\min(\bar\Gamma^{\max},\bar\Gamma_Q(\mu^*)))$ is feasible for~\eqref{form:gamma}, yielding $\gamma_Q(z)\leq \Gamma$.
\end{proof}
Observe that the above proposition considers the variant of \DDID\ defined by a strict inequality. This is not an issue in the case of \DDID\ as the MILP formulation (DDID-CG) provided in~\cite{OmerPossRougier2024} implies that encoding length of $\Phi(Q)$ is polynomially bounded in the encoding length of the instance (without the threshold $z$), so considering the strict inequality is polynomially equivalent to considering its non-strict variant. We formalize this observation in the following remark.
\begin{remark}
\label{rem:leq}
Proposition~\ref{prop:Phiandgamma} addresses the variant of~\DDID\ that seeks $Q\in\Q$ such that $\Phi(Q)<z$. This variant is polynomially equivalent to \DDID. Indeed, the threshold $z$ does not appear in formulation~\cite[(DDID-CG)]{OmerPossRougier2024}, so Lemma~\ref{eq:MILPinNP} shows that
$\Phi^* := \min_{Q\in\mathcal{Q}} \Phi(Q)$ is a rational number of encoding length at most $L$, where $L$ is
polynomial in the input size and does not depend on $z$. Hence, writing $\Phi^*=a/b$ and $z = r/s$, $b<2^L$, so we can define $\epsilon := 1/(s2^L)$ such that either $\Phi^* = z$ or $|\Phi^* - z| \ge \epsilon$. Therefore, $\Phi^* < z \iff \Phi^* \le z - \epsilon$ and $\Phi^* \le z \iff \Phi^* < z + \epsilon$, where $\epsilon$ has polynomial encoding length.
\end{remark}

\xc{In what follows, we study the complexity of \DDID\ under different choices of feasible sets  $\Q$ and $\Y$. For ease of reading, whenever one or both feasible sets $\Q$ and $\Y$ are specified,
we denote the corresponding \DDID\ problem by \ddid{$\Q$}{$\Y$}.}

We illustrate next the reformulation of Proposition~\ref{prop:Phiandgamma} and Remark~\ref{rem:leq} on the special case of \DDID, named \xc{\ddid{$\Q$}{$\Y^{\ONESEL}$}}, in which we define $\Y=\Y^{\ONESEL}:=\set{Y\subseteq[n]}{|Y|=1}=\set{\{i\}}{i\in[n]}$. 
Plugging the definition of $\Y^\ONESEL$ into~\eqref{form:gamma} leads to
\begin{subequations}
\label{form:1sel}
\begin{align}
\gamma_Q(z)\;:=\; \min_{\mu,\bar \Gamma}\; &\bar\Gamma+ \sum_{i\in Q}\frac{\mu_i}{d_i} \\
\text{s.t.}\; &z\leq c_i + \mu_i, &&\forall i\in Q,\\
&z\leq c_i + \min(1,\bar \Gamma)d_i, &&\forall i\in \bar Q,\\
& 0\leq \mu \leq d,\\
& 0\leq \bar\Gamma \leq 1,
\end{align}
\end{subequations}
where we assume that, without loss of generality, $\bQ$ is not empty.
Observe that if $z>\min_{i\in[n]}c_i+d_i$, problem~\eqref{form:1sel} is unfeasible so $\gamma_Q(z)=+\infty$ and $\Phi(Q)< z$. Otherwise, we can set $\mu_i=\pluspart{z-c_i}\leq d_i$ and $\bar\Gamma = \max_{i\in\bar Q}\pluspart{\frac{z-c_i}{d_i}}\leq 1$, so
\begin{equation}
\label{eq:ref_gamma_Q(z)}
 \gamma_Q(z) = \max_{i\in\bar Q}\frac{\pluspart{z-c_i}}{d_i} + \sum_{i\in Q}\frac{\pluspart{z-c_i}}{d_i}.   
\end{equation}
\xc{Consider any combinatorial optimization problem $\Pi$ for which an instance is given by a cost vector $w\in\Z^n$ and a set of feasible solutions $\Q^{\Pi}$ over which the total cost of the selected elements is to be optimized, formally:
\begin{equation}
\label{eq:Pi-OPT}
    \max_{Q\in \Q^{\Pi}}\; \sum_{i\in Q} w_i
    \tag{$\Pi$-OPT}
\end{equation}
Given a threshold $W$, its decision version can be expressed as
\begin{equation}
\label{eq:Pi}
    \max_{Q\in \Q^{\Pi}}\; \sum_{i\in Q} w_i> W.
    \tag{$\Pi$}
\end{equation}
Denote by \ddid{$\Q^{\Pi}$}{$\Y$} the resulting \DDID\ problem under the feasible set $\Q^{\Pi}$ of $\Pi$. We next show that \ddid{$\Q^{\Pi}$}{$\Y^{\ONESEL}$} is essentially as easy as $\Pi$.



\begin{proposition}
\label{prop:DDID-Pi}
    Given problem $\Pi$-\textsc{OPT} and the resulting problem \textup{\ddid{$\Q^{\Pi}$}{$\Y^{\ONESEL}$}},
    \begin{enumerate}
        \item[i.] if $\Pi$ is NP-hard, then \textup{\ddid{$\Q^{\Pi}$}{$\Y^{\ONESEL}$}} is NP-hard.
        \item[ii.] if $\Pi$ is polynomial-time solvable, then \textup{\ddid{$\Q^{\Pi}$}{$\Y^{\ONESEL}$}} is polynomial-time solvable.
    \end{enumerate}
\end{proposition}
\begin{proof}
Given any set $\Q$, the reformulation~\eqref{eq:ref_gamma_Q(z)} implies that, $\exists Q\in\Q$ such that $\gamma_Q(z) > \Gamma$ if and only if
\begin{equation}
\label{eq:ref_hard_Q}
\max_{Q\in\Q,\;\ell\in[n]\setminus Q}\frac{\pluspart{z-c_\ell}}{d_\ell} + \sum_{i\in Q}\frac{\pluspart{z-c_i}}{d_i} > \Gamma.
\end{equation} 

\noindent\textbf{i (hardness).}\; Consider any NP-hard problem $\Pi$. Deciding whether
\begin{equation}
\label{eq:hard_Q}
    \max_{Q\in \Q^{\Pi}}\; \sum_{i\in Q} w_i >W
\end{equation}
is thus hard, where $w_i>0$ is the cost of each $i\in [n]$. We now construct an instance of the resulting \ddid{$\Q^{\Pi}$}{$\Y^{\ONESEL}$}: let $\frac{[z-c_i]^+}{d_i}= \frac{w_i}{n\cdot(1+\max_{i\in [n]} w_i)}$ for each $i\in [n]$, and $\frac{[z-c_{n+1}]^+}{d_{n+1}}= 1$. There is no constraint on the item $n+1$. Therefore, there is a \textit{yes} solution to~\eqref{eq:hard_Q} if and only if the left-hand side of~\eqref{eq:ref_hard_Q} is greater than $\frac{W}{n\cdot(1+\max_{i\in [n]} w_i)}+1$. 

\noindent\textbf{ii (tractability).}\; 
First note that if $\Pi$ is polynomial-time solvable, so is $\Pi$-\textsc{OPT}. Given an instance of \ddid{$\Q^{\Pi}$}{$\Y^{\ONESEL}$}, we construct $n$ instances $\Pi(\ell)$-\textsc{OPT} of $\Pi$-\textsc{OPT}: for each $\ell \in[n]$, $w_\ell=-M$, with $M:=1+\sum_{i\in[n]}\frac{\pluspart{z-c_i}}{d_i}$, $w_i=\frac{\pluspart{z-c_i}}{d_i}$ for all $i\in [n]\setminus\{\ell\}$. Furthermore, $v(\ell):=\max_{Q\in\Q^\Pi}\sum_{i\in Q}w_i$ and $Q(\ell)$ denotes an optimal solution. Our algorithm for \ddid{$\Q^{\Pi}$}{$\Y^{\ONESEL}$} works as follows: we solve the $n$ problems $\Pi(\ell)$-\textsc{OPT} and return the solution $(Q(\tilde\ell),\tilde\ell)$ maximizing $z(\ell):=v(\ell)+\frac{\pluspart{z-c_\ell}}{d_\ell}$. If $z(\tilde \ell)>\Gamma$, output \emph{yes}, together with certificate $(Q(\tilde \ell),\tilde \ell)$; otherwise, \emph{no}.

Denote by $z^*>0$ the left-hand side of~\eqref{eq:ref_hard_Q}, and by $(Q^*,\ell^*)$ a maximizer. We prove the validity of the algorithm by showing that $z(\tilde \ell)=z^*$. Consider $\ell^*$. We see that
$$
\max_{Q\in\Q^\Pi}\sum_{i\in Q}w_i=\max\left\{\max_{Q\in\Q^\Pi: \ell^*\in Q}\sum_{i\in Q}w_i,\max_{Q\in\Q^\Pi: \ell^*\notin Q}\sum_{i\in Q}w_i\right\}
\geq 
\max_{Q\in\Q^\Pi: \ell^*\notin Q}\sum_{i\in Q}w_i
= 
z^*-\frac{\pluspart{z-c_\ell^*}}{d_\ell^*},
$$
where the inequality holds because the set $\set{Q\in\Q^\Pi}{\ell^*\notin Q}\neq\emptyset$ and the cost of any solution in $\set{Q\in\Q^\Pi}{\ell^*\notin Q}$ is larger than the cost of any solution in $\set{Q\in\Q^\Pi}{\ell^*\in Q}$. Furthermore, for each $\ell\in[n]$, either $\ell\in Q(\ell)$, so $z(\ell)\leq0$, or $\ell\notin Q(\ell)$, in which case $(Q(\ell),\ell)$ is a feasible solution for the left-hand side of~\eqref{eq:ref_hard_Q}, of cost $z(\ell)$; in both cases, $z(\ell)\leq z^*$. Hence, $z(\tilde \ell)= z^*$, so we can solve \ddid{$\Q^{\Pi}$}{$\Y^{\ONESEL}$} in polynomial-time.



\end{proof}


Therefore, if $\Q$ is hard-to-optimize, the resulting \DDID\ is NP-hard even when $\Y=\Y^{\ONESEL}$.
Henceforth, we focus on easy-to-optimize sets $\Q$, namely, those for which problem $\Pi$ is easy to solve. Proposition~\ref{prop:DDID-Pi}-(ii) shows that, for such sets $\Q$, the resulting \DDID\ is easy when considering the 
simplest feasible set $\Y^\ONESEL=\set{\{i\}}{i\in[n]}$. }
Despite the polynomiality of this special case, \DDID\ is hard in general, even for easy-to-optimize sets $\Q$. We show this in the next section, by considering the special case based on $\Y^\TWOSEL:=\set{\{i,j\}}{i,j\in[n],\ i\neq j}$.

\section{Hardness of \DDID}
\label{sec:hardness}

\xc{In this section, we establish the hardness of \ddid{$\Q$}{$\Y$} even for easy-to-optimize sets $\Q$.}

\subsection{Reformulation of $\Phi(Q)$ for $2$-selection}

We consider in what follows a special case of $\DDID$, named \xc{\ddid{$\Q^{\SEL}$}{$\Y^{\TWOSEL}$}}, for which $\Q=\Q^{\SEL}:=\set{Q\subseteq[n]}{|Q|=q}$, and $\Y=\Y^{\TWOSEL}:=\set{Y\subseteq[n]}{|Y|=2}=\set{\{i,j\}}{i,j\in[n],\ i\neq j}$ is the feasibility set of the 2-selection problem. An instance of \xc{\ddid{$\Q^{\SEL}$}{$\Y^{\TWOSEL}$}} is given by $n,q,c,d,\Gamma,z$.

Then, the expression $\max_{\xi\in\Xi_{\bar \Gamma}}\sum_{i\in Y\cap\bar Q}\xi_i$ in the right-hand side of~\eqref{form:gamma:Y} can be reformulated depending on $|Y\cap\bar Q|$:
\begin{equation}
  \label{eq:YcapbarQ}
\max_{\xi\in\Xi_{\bar \Gamma}}\sum_{i\in Y\cap\bar Q}\xi_i =
\begin{cases}
0 & \text{if }Y\cap\bar Q=\emptyset\\
\min(1,\bGam)\,d_j & \text{if }Y\cap\bar Q=\{j\}\\
\min(1,\bGam)d_i+\min\bigl(1,\max(0,\bGam-1)\bigr)d_j & \text{if }Y\cap\bar Q=\{i,j\},d_i\geq d_j.
\end{cases}
\end{equation}
Introducing
\begin{equation*}
  \bQ_2 = \set{(i,j)}{i,j\in\bQ,\ i\neq j,\ d_i\geq d_j}
\end{equation*}
we can plug~\eqref{eq:YcapbarQ} in the right-hand side of~\eqref{form:gamma:Y} to obtain
\begin{subequations}
  \label{form:gamma-2sel}
\begin{align}
\gamma_Q(z)\;=\; \min_{\mu,\bar \Gamma}\; &\bar\Gamma+ \sum_{i\in Q}\frac{\mu_i}{d_i} \\
\text{s.t.}\; & z \leq c_i+c_j+\mu_i+\mu_j, && \forall i,j\in Q,\ i\neq j, \label{eq:1a}\\
  & z \leq c_i+c_j+\mu_i+\min(1,\bGam)\,d_j, && \forall i\in Q,\ j\in\bQ, \label{eq:1b}\\
  & z \leq c_i+c_j+\min(1,\bGam)d_i+\min\bigl(1,\max(0,\bGam-1)\bigr)d_j,
    && \forall (i,j)\in\bQ_2, \label{eq:1c}\\
  & 0\leq \mu\leq d, \label{eq:1d}\\
  & 0\leq \bGam \leq 2. \label{eq:1e}
\end{align}
\end{subequations}
\subsection{Study of the optimal solutions of $\gamma_Q(z)$}

To simplify the presentation of the reduction provided in the next subsection, we introduce below properties which, whenever satisfied, leads to simplification of the function $\gamma_Q(z)$. Specifically, we consider only specific instances that satisfy $|Q|\leq n-2$ for each $Q\in\Q$, and
\begin{subequations}
  \label{eq:assumption}
\begin{align}
  &\bGam_{\min}(Q,z):=\max_{(i,j)\in\bQ_2}\ \frac{z-(c_i+c_j)}{d_i}=1,\label{eq:assumption:Gamma}\\
  &c_n+d_n<\tfrac{z}{2}.\label{eq:assumption:cndn}
\end{align}
\end{subequations}
We show next how assuming~\eqref{eq:assumption} leads to a closed-form reformulation for $\gamma_Q(z)$ if $n\in\bar Q$, while still involving one optimization variable otherwise. First, we show that this assumption implies we can assume that $\bar\Gamma=1$ in any optimal solution to~\eqref{form:gamma-2sel}.
\begin{proposition}\label{prop:1}
Assuming $|Q|\leq n-2$ and~\eqref{eq:assumption:Gamma},
\begin{align}
  \gamma_Q(z)\;=\;1+\min_{\mu}\quad & \sum_{i\in Q}\frac{\mu_i}{d_i}\nonumber\\
  \text{s.t.}\quad & \eqref{eq:1a} \text{ and } \eqref{eq:1d}\nonumber\\
  & c_i+c_j+\mu_i+d_j\;\geq\;z,\qquad \forall i\in Q,\ \forall j\in\bQ.\tag{\ref{eq:1b}$'$}\label{eq:prime}
\end{align}
\end{proposition}

\begin{proof}
  Assumption $\bGam_{\min}(Q,z)=1$ implies that
  $\min_{(i,j)\in\bQ_2}\; c_i+c_j+d_i = z$. Then \eqref{eq:1c} is satisfied for every
  $\bGam\geq 1$. If $\bGam< 1$, then the right-hand side of~\eqref{eq:1c} rewrites $c_i+c_j+\bGam d_i$, and
  $$
  \min_{(i,j)\in\bQ_2}\;c_i+c_j+\bGam d_i<\min_{(i,j)\in\bQ_2}\; c_i+c_j+d_i = z,
  $$
  so~\eqref{eq:1c} is violated for some $(i,j)\in\bQ_2$.
 
Rewriting \eqref{eq:1b} for $\bGam\geq 1$ yields \eqref{eq:prime}, so variable $\bar\Gamma$ only appears in the objective and its bound constraints. It follows that the minimum is
attained at $\bGam=1$.
\end{proof}

The above reformulation is leveraged in the next result to provide an expression for $\gamma_Q(z)$ that involves at most one optimization variable. 

\begin{proposition}\label{prop:2}
Let $Q\in\Q$ be such that $|Q|\leq n-2$,~\eqref{eq:assumption} holds and~\eqref{form:gamma-2sel} is feasible. Then:
\begin{description}
\item[(a)] If $n\in\bar Q$:
  \begin{equation*}
    \gamma_Q(z)=1+\sum_{i\in Q}\frac{\pluspart{z-(c_n+d_n+c_i)}}{d_i}.
  \end{equation*}
\item[(b)] If $n\in Q$:
  \begin{align}
    \gamma_Q(z)=1+\min_{\mu_n}\quad
      & \frac{\mu_n}{d_n}+\sum_{i\in Q\setminus\{n\}}
        \frac{\pluspart{z-(c_n+\mu_n+c_i)}}{d_i}\\
    \text{s.t.}\quad
      & \mu_n \geq\; z-c_n-\min_{j\in[n-1]}\,(c_j+d_j)\label{eq:lowerboundmun}\\
      & 0\leq \mu_n\leq d_n .
  \end{align}
\end{description}
\end{proposition}

\begin{proof}
We first note that~\eqref{eq:assumption:cndn} and the feasibility of~\eqref{form:gamma-2sel} imply that
\begin{subequations}
  \label{eq:star}
\begin{align}
    c_i+d_i>\frac{z}{2}, && \forall i\in\bQ\setminus\{n\},\label{eq:star:d}\\
    c_i+\mu_i>\frac{z}{2}, && \forall i\in Q\setminus\{n\},\label{eq:star:mu}
\end{align}
\end{subequations}
for any feasible solution $(\mu,\bar\Gamma)$ to~\eqref{form:gamma-2sel}. Indeed, assume there is $i\in\bQ \setminus\{n\}$ such that $c_i+d_i\leq \frac{z}{2}$. Then, if $n\in \bQ$, $d_n\leq d_i$ (the case $d_i\leq d_n$ is proved similarly), we have that
$$
c_i+c_n+\min(1,\bGam)d_i+\min\bigl(1,\max(0,\bGam-1)\bigr)d_n\leq c_i+c_n+d_i+d_n < z,
$$
so~\eqref{eq:1c} is infeasible. If $n \in Q$, relying instead on~\eqref{eq:1b} yields the same conclusion. If $i\in Q\setminus\{n\}$, \eqref{eq:star:mu} is proved by relying instead on~\eqref{eq:1a} and~\eqref{eq:1b}. 

In what follows, we let $\mu^{*}$ be an optimal solution.

\medskip
\noindent\textbf{(a) Case $n\in \bar Q$.}
 By~\eqref{eq:star:mu},
\begin{equation*}
  c_i+c_j+\mu^{*}_i+\mu^{*}_j>z,\qquad \forall i,j\in Q,\ i\neq j,
\end{equation*}
so that \eqref{eq:1a} is satisfied strictly. From~\eqref{eq:star}, we also have
\begin{equation*}
  c_i+c_j+\mu^{*}_i+d_j>z,\qquad \forall i\in Q,\ \forall j\in\bQ\setminus\{n\},
\end{equation*}
so~\eqref{eq:prime} is also satisfied strictly.
Apart from the bounds on $\mu^{*}$, the only remaining constraints are
\begin{equation*}
  c_i+c_n+\mu^{*}_i+d_n \geq z .
\end{equation*}
By optimality of $\mu^{*}$, we see that for each $i\in Q$, either $c_i+c_n+\mu^{*}_i+d_n = z$ or $\mu^{*}_i=0$, for otherwise decreasing $\mu^{*}_i$ by some $\epsilon>0$ yields a strictly better solution. Hence, for
every $i\in Q$,
$\mu^{*}_i=\pluspart{z-(c_i+c_n+d_n)}$, so that
\begin{equation*}
  \gamma_Q(z)\;=\;1+\sum_{i\in Q}\frac{\pluspart{z-(c_n+d_n+c_i)}}{d_i}.
\end{equation*}

\medskip
\noindent\textbf{(b) Case $n\in Q$.} Similarly to (a), one shows that
the constraints not involving item $n$ are not active so the only constraint (other than bounds) involving $\mu^*_i$ that may possibly be active is~\eqref{eq:1a}, 
\begin{equation*}
  c_n+c_i+\mu^*_n+\mu^*_i \geq z.
\end{equation*}
As before, for every $i\in Q\setminus\{n\}$, we obtain
\begin{equation*}
  \mu^{*}_i=\pluspart{z-(c_n+\mu^{*}_n+c_i)}\leq d_i.
\end{equation*}
The constraints \eqref{eq:prime} involving item $n$ give
\begin{align*}
  & c_n+\mu^{*}_n+c_j+d_j\;\geq\;z, && \forall j\in\bQ\\
  \iff\ & \mu^{*}_n\;\geq\;z-c_n-(c_j+d_j), && \forall j\in\bQ\\
  \iff\ & \mu^{*}_n\;\geq\;z-c_n-\min_{j\in\bQ}(c_j+d_j). &&
\end{align*}
To conclude, we note that $\mu^{*}_i\leq d_i$ for every $i\in Q\setminus\{n\}$
implies
\begin{align*}
  & z-(c_n+\mu^{*}_n+c_i)\;\leq\;d_i, && \forall i\in Q\setminus\{n\}\\
  \iff\ & \mu^{*}_n\;\geq\;z-c_n-(c_i+d_i), && \forall i\in Q\setminus\{n\}.
\end{align*}
We deduce that
\begin{equation*}
\begin{aligned}
  \gamma_Q(z)\;=\;1+\min_{\mu_n}\quad
    & \frac{\mu_n}{d_n}+\sum_{i\in Q\setminus\{n\}}
      \frac{\pluspart{z-(c_n+\mu_n+c_i)}}{d_i}\\
  \text{s.t.}\quad
    & \mu_n\;\geq\;z-c_n-\min_{i\in[n-1]}(c_i+d_i)\\
    & 0\leq\mu_n\leq d_n.
\end{aligned}
\end{equation*}
\end{proof}

\subsection{Polynomial reduction}

The hardness proof is carried out by a polynomial reduction from
\PARTITION, which, given a list of integers $(h_1,\dots,h_{2k})\subseteq\N$ summing up to $2H$ for some integer $H$, asks whether there exist $k$ elements with indices $T\subseteq[2k]$ such that
\begin{equation*}
  \sum_{i\in T} h_i = H.
\end{equation*}
\PARTITION\ is known to be NP-hard~\cite[SP12]{GareyJohnson1979}. We consider without loss of generality that $k\geq 12$, and reduce any instance $I$ of \PARTITION\ to an instance $I'$ of \xc{\ddid{$\Q^{\SEL}$}{$\Y^{\TWOSEL}$}}. Following Proposition~\ref{prop:Phiandgamma} and Remark~\ref{rem:leq}, we consider the variant of \xc{\ddid{$\Q^{\SEL}$}{$\Y^{\TWOSEL}$}} that seeks $Q\in\Q^\SEL$ such that $\gamma_Q(z)>\Gamma$. 

The instance $I'$ contains $n:=3k+4$ items and $q:=k+1$. Each of the first $2k$ items $i\in[2k]$ corresponds to an element of \PARTITION,  while the last $k+4$ items are added so as to satisfy the assumptions of Proposition~\ref{prop:2}, namely \eqref{eq:assumption} and ensuring the feasibility of~\eqref{form:gamma-2sel}. Specifically, letting $\omega_i=h_i+2H$ for each $i\in[2k]$, and $\Omega=(2k+1)H$, we further define
\begin{align*}
  \forall i\in[2k]: \quad
    & c_i = \frac{z}{2}+q\cdot d_n-\Bigl(\tfrac{3}{4}d_i+\tfrac{1}{2}d_n\Bigr),
    && d_i=\frac{\Omega}{\omega_i}\,d_n \;>\; d_n,\\
  \forall g\in G:=\{2k+1,\dots,2k+q+2\}: \quad
    & c_g = 0, && d_g = z,\\
  & c_n = \frac{z}{2}-q\cdot d_n, && d_n = \sumh,\\
  & z = 6(q+1)\,d_n, &&
\end{align*}
and $\Gamma:=\frac{3}{2}+\frac{3}{4}(q-1)-\frac{1}{2\Omega}$. Observe that for each $i\in[2k]$ and $k\geq 12$ the constructed instance satisfies
\begin{equation}
d_i\geq 6d_n\label{digeqdn}.
\end{equation}
We show next that the constructed instance satisfies all requirements of Proposition~\ref{prop:2}. Property~\eqref{eq:assumption:cndn} and $|Q|\leq n-2$ for each $Q\in\Q$ follow immediately from the definitions of $q$, $c_n$, and $d_n$, and we prove below that~\eqref{eq:assumption:Gamma} holds for all $Q\in\Q^\SEL$. 
\begin{proposition}\label{prop:3}
For every $Q\in\Q^\SEL$, we have $\bGam_{\min}(Q,z)=1$.
\end{proposition}
\begin{proof}
Let $Q\in\Q^\SEL$. Since $|G|=|Q|+2$, there
exist $g_1,g_2\in G$, $g_1\neq g_2$, with $g_1\in\bQ$ and $g_2\in\bQ$. To compute $\bGam_{\min}(Q,z)$ we evaluate $\frac{z-(c_i+c_j)}{d_i}$ for all
$(i,j)\in\bQ_2$, enumerating the cases. We start with the case where $i,j\in G$, for which $\frac{z-(c_i+c_j)}{d_i}=\frac{z}{z}=1$.

To compute the remaining cases, we first show that for each $i\in[2k]$,
\begin{equation}
c_i\geq \frac{z}{2}\label{eq:cigeq0}
\end{equation}
follows from (i) $\tfrac{3}{4}d_i+\tfrac{1}{2}d_n\leq \tfrac{3}{4}d_i+\tfrac{1}{12}d_i\leq d_i$ (using~\eqref{digeqdn}), and (ii) $q\cdot d_n=(k+1)\,\frac{\omega_i}{\Omega}d_i=(k+1)\frac{h_i+2H}{(2k+1)H}d_i\geq\frac{(2k+1)H}{(2k+1)H}d_i$. We now consider the three different cases:
\begin{itemize}
  \item $i\in G$, $j\in[2k]\cup\{n\}$:
    $\displaystyle \frac{z-(c_i+c_j)}{d_i}\leq\frac{z}{d_g}\leq 1$, from~\eqref{eq:cigeq0} if $j\in[2k]$, and using $d_g>d_i>d_n$ for each $i\in[2k]$,
  \item $i\in[2k]$, $j=n$:
    $\displaystyle
      \frac{z-(c_i+c_j)}{d_i}
      =\frac{3}{4}+\frac{1}{2}\,\frac{d_n}{d_i}
      \leq 1,
    $
    from~\eqref{digeqdn},
  \item $i,j\in[2k]$: $\frac{z-(c_i+c_j)}{d_i}
      \leq 0$, from~\eqref{eq:cigeq0}.
    \qedhere
\end{itemize}
\end{proof}
We prove next that the feasibility of~\eqref{form:gamma-2sel} for all $Q\in\Q^\SEL$ is a direct consequence of Proposition~\ref{prop:3}.
\begin{proposition}\label{prop:new}
  For each $Q\in\Q^\SEL$, problem~\eqref{form:gamma-2sel} is feasible.
\end{proposition}
\begin{proof}
We set $\bar \Gamma=2$ and $\mu_i=d_i$ for all $i\in Q$, so the problem is feasible if $c_i+c_j+d_i+d_j\geq z$ for each $i,j\in[n]$, $i\neq j$. This result follows immediately from $\frac{z-(c_i+c_j)}{d_i}\leq 1$ for each $i,j\in[n]$, $d_i\geq d_j$, proved in Proposition~\ref{prop:3}.
\end{proof}

In the following proposition, we show that it is suboptimal to observe the elements of $G$.

\begin{proposition}\label{prop:4}
For every $Q_G\in\Q$ there exists $Q\in\Q$ with $Q\subseteq[2k]\cup\{n\}$ such that
$
  \gamma_Q(z)\;\geq\;\gamma_{Q_G}(z).
$
\end{proposition}

\begin{proof}
If $Q_G\cap G=\emptyset$, it suffices to take $Q=Q_G$.

Suppose now $Q_G\cap G\neq\emptyset$ and let $g\in Q_G\cap G$. For
$0\leq\mu_n\leq d_n$,
\begin{equation}
  \label{eq:ub23}
  \frac{\pluspart{z-(c_n+\mu_n+c_g)}}{d_g}
  =\frac{\pluspart{z-(c_n+\mu_n)}}{z}
  =\frac{1}{2}+\frac{q\cdot d_n-\mu_n}{z}\leq \frac{1}{2}+\frac{q\cdot d_n}{6(q+1)d_n}
  \leq\frac{1}{2}+\frac{1}{6}=\frac{2}{3}.
\end{equation}
For $j\in[2k]\setminus Q_G$, by~\eqref{digeqdn} we have
\begin{equation}
  \label{eq:lb23}
  \frac{\pluspart{z-(c_n+\mu_n+c_j)}}{d_j}
  =\frac{\pluspart{\frac{3}{4}d_j+\frac{1}{2}d_n-\mu_n}}{d_j}
  \geq\frac{3}{4}-\frac{1}{2}\,\frac{d_n}{d_j}
  \geq\frac{2}{3}.
\end{equation}

Let $I\subseteq[2k]\setminus Q_G$ with $|I|=|Q_G\cap G|$ (such a set exists for
$2k\geq q=k+1$), and set $Q=(Q_G\setminus G)\cup I$. Referring to Proposition~\ref{prop:2}, we have two cases:

\medskip
\noindent\textbf{(a) Case $n\in \bar Q_G$.}
\begin{align*}
  \gamma_{Q_G}(z)
  &=1+\sum_{g\in Q_G\cap G}\frac{\pluspart{z-(c_n+d_n+c_g)}}{d_g}
     +\sum_{i\in Q_G\setminus G}\frac{\pluspart{z-(c_n+d_n+c_i)}}{d_i}\\
  &\leq 1+\frac{2}{3}\,\bigl|Q_G\cap G\bigr|
     +\sum_{i\in Q_G\setminus G}\frac{\pluspart{z-(c_n+d_n+c_i)}}{d_i}.
\end{align*}
 Furthermore,
$\displaystyle
  \gamma_{Q}(z)\;\geq\;1+\frac{2}{3}\,|I|
    +\sum_{i\in Q_G\setminus G}\frac{\pluspart{z-(c_n+d_n+c_i)}}{d_i}
    $. Hence, $\gamma_Q(z)\;\geq\;\gamma_{Q_G}(z)$.

\medskip
\noindent\textbf{(b) Case $n\in Q_G$.} Bounds~\eqref{eq:ub23} and~\eqref{eq:lb23} imply that
\begin{equation*}
  \sum_{i\in Q\setminus\{n\}}\frac{\pluspart{z-(c_n+\mu_n+c_i)}}{d_i}
  \;\geq\;
  \sum_{i\in Q_G\setminus\{n\}}\frac{\pluspart{z-(c_n+\mu_n+c_i)}}{d_i},
\end{equation*}
for each $0\leq \mu_n\leq d_n$. Taking the minimum over $\mu_n$ leads to
$\gamma_Q(z)\geq\gamma_{Q_G}(z)$.
\end{proof}
We now focus on sets $Q$ that do not contain elements of $G$, and compute a closed-form expression for $\gamma_Q$.
\begin{proposition}\label{prop:5}
Let $Q\subseteq[2k]\cup\{n\}$ with $|Q|=q$.
\begin{itemize}
  \item If $n\in \bar Q$:
    $\displaystyle
      \gamma_Q(z)=1+\frac{3}{4}q-\frac{1}{2}\sum_{i\in Q}\frac{\omega_i}{\Omega}. $
  \item If $n\in Q$:
    $\displaystyle
      \gamma_Q(z)=\frac{3}{2}+\frac{3}{4}(q-1)
      -\frac{1}{2}\left|1-\sum_{i\in Q\setminus\{n\}}\frac{\omega_i}{\Omega}\right| . $
\end{itemize}
\end{proposition}

\begin{proof}
We start from Proposition~\ref{prop:2}.

\medskip
\noindent\textbf{(a) Case $n\in \bar Q$.}
\begin{equation*}
  \gamma_Q(z)
  =1+\sum_{i\in Q}\frac{\frac{3}{4}d_i-\frac{1}{2}d_n}{d_i}
  =1+\sum_{i\in Q}\frac{\bigl(\frac{3}{4}-\frac{1}{2}\frac{\omega_i}{\Omega}\bigr)d_i}{d_i}
  =1+\frac{3}{4}q-\frac{1}{2}\sum_{i\in Q}\frac{\omega_i}{\Omega}.
\end{equation*}

\medskip
\noindent\textbf{(b) Case $n\in Q$:} We first show that with chosen parameters, constraints \eqref{eq:lowerboundmun} are redundant in the definition of $\gamma_Q$. Observe that
\begin{align*}
  \min_{j\in[n-1]}(c_j+d_j)
  &=\min_{j\in[2k]}\ \frac{z}{2}+q\cdot d_n-\frac{3}{4}d_j-\frac{1}{2}d_n+d_j\\
  &=\frac{z}{2}+\Bigl(q-\frac{1}{2}\Bigr)d_n+\frac{1}{4}\min_{j\in[2k]} d_j.
\end{align*}
Using~\eqref{digeqdn},
\begin{align*}
  z-c_n-\min_{j\in[n-1]}(c_j+d_j)
  &= z-\frac{z}{2}+q\cdot d_n-\frac{z}{2}-q\cdot d_n+\frac{1}{2}d_n
     -\frac{1}{4}\min_{j\in[2k]}d_j\\
  &= \frac{1}{2}d_n-\frac{1}{4}\min_{j\in[2k]}d_j
  \leq0.
\end{align*}
Therefore, the lower bound~\eqref{eq:lowerboundmun} are dominated by the non-negativity constraint on $\mu_n$, so
\begin{align*}
  \gamma_Q(z)&=1+\min_{0\leq\mu_n\leq d_n}\ \frac{\mu_n}{d_n}
  +\sum_{i\in Q\setminus\{n\}}
    \left(\frac{3}{4}+\frac{1}{2}\frac{d_n}{d_i}-\frac{\mu_n}{d_i}\right),\\
  &=1+\frac{3}{4}(q-1)+\sum_{i\in Q\setminus\{n\}}\frac{1}{2}\frac{d_n}{d_i}
    +\min_{0\leq\mu_n\leq d_n}\ \frac{\mu_n}{d_n}
      \left(1-\sum_{i\in Q\setminus\{n\}}\frac{\omega_i}{\Omega}\right)\\[1ex]
  &=1+\frac{3}{4}(q-1)+\frac{1}{2}\sum_{i\in Q\setminus\{n\}}\frac{\omega_i}{\Omega}
    +\begin{cases}
       1-\sum_{i}\frac{\omega_i}{\Omega}, & \text{if } \sum_{i\in Q\setminus\{n\}} \omega_i>\Omega,\\
       0, & \text{otherwise,}
     \end{cases}\\[1ex]
  &=\frac{3}{2}+\frac{3}{4}(q-1)
    -\frac{1}{2}\left|1-\sum_{i\in Q\setminus\{n\}}\frac{\omega_i}{\Omega}\right|.
\end{align*}
\end{proof}

\begin{theorem}\label{thm:main}
$(h_1,\dots,h_{2k})$ is a yes-instance of \PARTITION{} if and only if $\exists Q\in\Q^\SEL$ such that 
$
  \gamma_Q(z)
  >\frac{3}{2}+\frac{3}{4}(q-1)-\frac{1}{2\Omega}.
$
\end{theorem}

\begin{proof}
Proposition~\ref{prop:4} shows that the search space for $Q$ can be restricted
to $[2k]\cup\{n\}$. By Proposition~\ref{prop:5}, and
noting that $\sum_{i\in Q}\omega_i\geq(k+1)2H\geq\Omega$, 
we get for $n\in \bar Q$ that
\begin{equation*}
  \gamma_Q(z)\leq\frac{1}{2}+\frac{3}{4}q
  =\frac{5}{4}+\frac{3}{4}(q-1)
  \leq\frac{3}{2}+\frac{3}{4}(q-1)-\frac{1}{2\Omega},
\end{equation*}
since $\Omega \geq 2$.

For $n\in Q$ we have that
$\gamma_Q(z)=\frac{3}{2}+\frac{3}{4}(q-1)$ if and only if
$\sum_{i\in Q\setminus\{n\}}\omega_i=\Omega$. In particular, if $\sum_{i\in Q\setminus\{n\}}\omega_i\neq\Omega$, we have from the definition of $\omega_i$ and $\Omega$ that
\begin{equation*}
  \left|1-\sum_{i\in Q\setminus\{n\}}\frac{\omega_i}{\Omega}\right|
  =\left|\frac{H-\sum_{i\in Q\setminus\{n\}} h_i}{\Omega}\right|
  \geq\frac{1}{\Omega},
\end{equation*}
so that
\begin{equation*}
  \gamma_Q(z)\;\leq\;\frac{3}{2}+\frac{3}{4}(q-1)-\frac{1}{2\Omega}.
\end{equation*}

Therefore, there exists $Q\in\Q^\SEL$ such that 
$\gamma_Q(z)
  >\frac{3}{2}+\frac{3}{4}(q-1)-\frac{1}{2\Omega}
  $
  if and only if there exists $I\subseteq[2k]$ with $|I|=k$ such that $\sum_{i\in I} \omega_i=\Omega$, which is equivalent to $\sum_{i\in I}h_i = H$ since $\sum_{i\in I} \omega_i=2k\cdot H+\sum_{i\in I} h_i$ and $\Omega=(2k+1)H$. 
\end{proof}

\begin{corollary}
    \textup{\ddid{$\Q^{\SEL}$}{$\Y^{\TWOSEL}$} is NP-complete.}
\end{corollary}
\begin{proof}
The hardness follows immediately from the Proposition~\ref{prop:Phiandgamma}, Remark~\ref{rem:leq} and Theorem~\ref{thm:main}, as well as the hardness of \PARTITION. The membership is proved in Corollary~\ref{cor:NP}.
\end{proof}

\xc{Hence, we conclude that \DDID\ is NP-hard even when $\Q=\Q^{\SEL}$ and $|\Y|$ is only $O(n^2)$.}

\section{\ddid{$\Q^{\SEL}$}{$\Y^{\SEL}$} with equal deviations}
\label{sec:equaldev}
In this section, we focus on the optimization variant of \xc{\ddid{$\Q^{\SEL}$}{$\Y^{\SEL}$}, where $\Q^{\SEL}:=\set{Q\subseteq[n]}{|Q|=q}$, and $\Y^{\SEL}:=\set{Y\subseteq[n]}{|Y|=p}$, under the assumption that $d_i=d_j$ for all $i,j\in [n]$.} We denote this problem by \xc{\ddid{$\Q^{\SEL}$}{$\Y^{\SEL}$}-\Equal-\textup{\textsf{OPT}}} and show that it can be solved in polynomial time. To establish its polynomial-time solvability, we begin with a simplified reformulation for this case, in which the variables $\xi_i$ for all $i\in \bQ$ in the innermost maximization are eliminated. Two observations are then stated that allow us to easily identify a worst-case scenario $\mu$ for any given $Q\in \Q^{\SEL}$, as well as an optimal solution $y$ for any given $Q\in \Q^{\SEL}$ and $\mu\in \Xi_{\Gamma}$. Building on these preliminary results, we provide two lemmas to reduce the set of candidate optimal solution $Q$ to polynomial size, thereby establishing polynomial-time solvability for this case. 

For ease of expressions, we assume that the elements are indexed in non-decreasing order of $c$, i.e., $c_1\leq c_2\leq\ldots\leq c_n$. Denote $d_i=d$ for all $i\in [n]$. 
Given $Q\in \Q^{\SEL}$ and $\mu\in \Xi_{\Gamma}$ we define function
$$
\varphi(Q,\mu):=\min_{Y\in \Y^\SEL}
         \max_{\xi \in \Xi_\Gamma(Q,\obsxi)} \sum_{i\in Y}(c_i+\xi_i).
$$
We first introduce the following reformulation of $\varphi(Q,\mu)$ for this special case:

$$
\varphi'(Q,\mu):=\min_{Y\in \Y^\SEL}
         \sum_{i\in Y\cap Q}(c_i+\mu_i)+ \sum_{i\in Y\cap \bar Q}c_i+ d\cdot\min\left\{|Y\cap \bar Q|, \; \bar{\Gamma}_Q(\mu)\right\}.
$$

\begin{lemma}
\label{lemma:ref_d_equal}
    For each $Q\in \Q^{\SEL}$ and $\mu\in \Xi_{\Gamma}$, we have $\varphi'(Q,\mu)=\varphi(Q,\mu)$.
\end{lemma}
\begin{proof}
    For each $Q\in \Q^{\SEL}$ and $\mu\in \Xi_{\Gamma}$, we have
    \begin{align*}
    \varphi(Q,\mu)&=\min_{Y\in \Y^\SEL}
         \max_{\xi \in \Xi_\Gamma(Q,\obsxi)} \sum_{i\in Y}(c_i+\xi_i)\\
         &=\min_{Y\in \Y^\SEL}
         \max_{\xi \in \Xi_\Gamma(Q,\obsxi)} \sum_{i\in Y\cap Q}(c_i+\mu_i)+\sum_{i\in Y\cap \bar Q}(c_i+\xi_i).
    \end{align*}
    To prove $\varphi'(Q,\mu)=\varphi(Q,\mu)$, it suffices to show that, for any given $Y\in \Y^\SEL$, $\max_{\xi \in \Xi_\Gamma(Q,\obsxi)}\sum_{i\in Y\cap \bar Q}\xi_i=d\cdot\min\left\{|Y\cap \bar Q|, \; \bar{\Gamma}_Q(\mu)\right\}$. Since $\xi_i\leq d$ for all $i\in Y\cap \bar Q$, we have $\sum_{i\in Y\cap \bar Q}\xi_i\leq d\cdot|Y\cap \bar Q|$. Since $\sum_{i\in \bar Q}\xi_i\leq d\cdot\bar{\Gamma}_Q(\mu)$, we also have $\sum_{i\in Y\cap \bar Q}\xi_i\leq  d\cdot\bar{\Gamma}_Q(\mu)$. Thus, $\max_{\xi \in \Xi_\Gamma(Q,\obsxi)}\sum_{i\in Y\cap \bar Q}\xi_i\leq d\cdot\min\left\{|Y\cap \bar Q|, \; \bar{\Gamma}_Q(\mu)\right\}$.
    We next show the reverse inequality. When $|Y\cap \bar Q|\leq \bar{\Gamma}_Q(\mu)$,
    let $\xi^*_i=d$ for $i\in Y\cap \bar Q$ and $\xi^*_k=0$ for $k\in \bar Q\setminus Y$. Then, $\sum_{i\in Y\cap \bar Q}\xi^*_i=d\cdot |Y\cap \bar Q|$. When $|Y\cap \bar Q|> \bar{\Gamma}_Q(\mu)$, let $\xi^*_i=d$ for $i\in S$, where $S\subset (Y\cap \bar Q)$ with $|S|=\lfloor \bar{\Gamma}_Q(\mu) \rfloor$,  $\xi^*_j=d\cdot(\bar{\Gamma}_Q(\mu)- \lfloor \bar{\Gamma}_Q(\mu) \rfloor)$ for some $j\in (Y\cap \bar Q)\setminus S$, and $\xi^*_k=0$ for $k\in \bQ \setminus (S\cup \{j\})$. Then, $\sum_{i\in Y\cap \bar Q}\xi^*_i=d \cdot\bar{\Gamma}_Q(\mu)$. Thus, $\max_{\xi \in \Xi_\Gamma(Q,\obsxi)}\sum_{i\in Y\cap \bar Q}\xi_i\geq d\cdot\min\left\{|Y\cap \bar Q|, \; \bar{\Gamma}_Q(\mu)\right\}$.
\end{proof}

Lemma~\ref{lemma:ref_d_equal} yields the following reformulation of \xc{\ddid{$\Q^{\SEL}$}{$\Y^{\SEL}$}-\Equal-\textup{\textsf{OPT}}}:
\begin{equation}
\label{eq:re_ddid_equal_dev}
\min_{Q\in \Q^\SEL} \max_{\obsxi\in \Xi_\Gamma} \min_{Y\in \Y^\SEL}
         \sum_{i\in Y\cap Q}(c_i+\mu_i)+ \sum_{i\in Y\cap \bar Q}c_i+ d\cdot\min\left\{|Y\cap \bar Q|, \; \bar{\Gamma}_Q(\mu)\right\}.
\end{equation}

For ease of notation, we henceforth write $\Q$ and $\Y$ instead of $\Q^\SEL$ and $\Y^\SEL$, respectively.

We next state two observations, the first one characterizes a worst-case scenario $\mu$, while the second one characterizes an optimal solution $Y$ for the components in $\bQ$.

\begin{observation}
\label{obs:opt_Q}
    Given $Q\in\Q$, there always exists a scenario $\mu^* \in \argmax_{\mu \in \Xi_{\Gamma}}\varphi(Q,\mu)$ such that $c_i+\mu^*_i\leq c_j+\mu^*_j$ for all $i,j\in Q$ with $i< j$.
\end{observation}
\begin{proof}
    Define $\mu'\in \argmax_{\mu \in \Xi_{\Gamma}}\varphi(Q,\mu)$ satisfying that $c_i+\mu'_i> c_j+\mu'_j$. We now construct a scenario $\mu^*$ through 
    $$
    \mu^*_j:=c_i-c_j+\mu'_i,
    $$ 
    $$\mu^*_i:=\mu'_i+\mu'_j-\mu^*_j \text{ and } \mu^*_k:=\mu'_k \text{ for all } k\neq i,j.$$
    Since $0\leq\mu'_j<\mu^*_j\leq\mu'_i\leq d$, $0\leq\mu'_j\leq\mu^*_i<\mu'_i\leq d$, and $\sum_{k\in[n]}\mu^*_k=\sum_{k\in[n]}\mu'_k$, the constructed scenario $\mu^*\in \Xi_{\Gamma}$. Then, we have 
    \begin{align*}
        c_i+ \mu^*_i&=c_i+(\mu'_i+\mu'_j-\mu^*_j)\\
        &=c_i+(\mu'_i+\mu'_j-c_i+c_j-\mu'_i)\\
        &=c_j+\mu'_j.
    \end{align*}
    Therefore, we have $c_i+ \mu^*_i=c_j+\mu'_j<c_i+\mu'_i=c_j+\mu^*_j$.
    Noticing that $\bar \Gamma_Q(\mu^*)=\bar \Gamma_Q(\mu')$, we conclude that $\varphi(Q,\mu^*)=\varphi(Q,\mu')$.
\end{proof}

\begin{observation}
\label{obs:opt_barQ}
    Given $Q\in \Q$ and $\mu\in \Xi_{\Gamma}$, there is an optimal solution $Y^*$ to the minimization problem defining $\varphi(Q,\mu)$ such that, 
    for all $i,j\in \bar Q$ with $i< j$, $j\in Y^*$ implies $i\in Y^*$.
\end{observation}
\begin{proof}
    Define $Y'\in \Y$ by setting $i\notin Y'$, $j\in Y'$ and $Y'\setminus \{i,j\}=Y^*\setminus \{i,j\}$,
    where $i,j\in \bar Q$ with $i< j$. The difference in costs of the solutions $Y'$ and $Y^*$ is
    \begin{align*}
        &\left(\sum_{k\in Y'\setminus\{i,j\}}c_k+c_j+\sum_{r\in Y'\cap Q} \mu_r+d\cdot\min\left\{|Y'\cap \bQ|, \; \bar{\Gamma}_Q(\mu)\right\}\right)\\
        &\quad \quad -    \left(\sum_{k\in Y^*\setminus\{i,j\}}c_k+c_i+\sum_{r\in Y^*\cap Q}\mu_r+d\cdot\min\left\{|Y^*\cap \bQ|,\; \bar{\Gamma}_Q(\mu)\right\}\right)=c_j-c_i
        \geq 0.
    \end{align*}
\end{proof}

\begin{remark}
\label{remark:opt_y}
    Given $Q^*\in\Q$ and $\mu^*\in\argmax_{\mu\in \Xi_{\Gamma}}\varphi(Q^*,\mu)$, there always exists an optimal solution $Y^*$ to the minimization problem defining $\varphi(Q^*,\mu^*)$ such that $Y^*$ contains the $s$ elements with the smallest values of $c+\mu$ in $Q^*$ and the first $(p-s)$ elements in $\bQ^*$, where $s\in [\min\{p,q\}]_0$.
\end{remark}

Based on the preliminary results stated above, we now establish two essential results, which allow us to eliminate suboptimal solutions $Q$, leading to the polynomial-time solvability for this special case. Specifically, among the elements in $[n]\setminus[p]$, any solution in which an element of $\bQ$ precedes an element of $Q$ is suboptimal; whereas the reverse ordering holds for the elements in $[p]$. 

The main idea behind these two proofs is to exchange the first element that violates the prescribed ordering with the closest element that restores the required ordering. We then construct a scenario $\mu'\in \Xi_{\Gamma}$ for the suboptimal solution $Q'$ based on a worst-case scenario $\mu^*$ for the solution obtained after this exchange, denoted by $Q^*$. Finally, we show that, for every candidate optimal solution $Y'$ under $Q'$ and $\mu'$, there exists a feasible solution $Y$ under $Q^*$ and $\mu^*$ whose objective value is no greater than that of $Y'$, thereby proving the suboptimality of $Q'$. 

\begin{lemma}
\label{lemma:opt_bar_p}
    There is an optimal solution $Q^*\in \Q$ such that $i< j$ for all $i\in Q^*\setminus[p]$ and $j\in \bQ^*\setminus[p]$.
\end{lemma}
\begin{proof}
    Denote 
    $$
    Q'=\{q_1,\ldots, q_2,j,q_3,\ldots,q_4\},
    $$
    and 
    $$
    \bQ'=\{\bar{q}_1,\ldots, \bar{q}_2,i,\ldots,j-1,\bar{q}_3\ldots,\bar{q}_4\},
    $$
    where $i,j\in [n]\setminus[p]$, $i$ is the first element in $\bQ'\setminus[p]$, and $j$ is the element in $Q'$ with the smallest index greater than $i$. Define
    $$
    Q^*=Q'\setminus\{j\}\cup \{i\}, \text{ and } \bQ^*=\bQ'\setminus\{i\}\cup \{j\},
    $$
    namely
    $$
    Q^*=\{q_1,\ldots, q_2,i,q_3,\ldots,q_4\},
    $$
    and 
    $$
    \bQ^*=\{\bar{q}_1,\ldots, \bar{q}_2,\ldots,j-1,j,\bar{q}_3\ldots,\bar{q}_4\}.
    $$
    Let $\mu^*\in \argmax_{\mu \in \Xi_{\Gamma}}\varphi(Q^*,\mu)$ satisfying Observation~\ref{obs:opt_Q}, and, without loss of generality, assume that $\mu^*_k=0$ for each $k\in\bar Q^*$.
    We next show that $\Phi(Q^*)\leq \Phi(Q')$. To do so, we construct a scenario $\mu'\in \Xi_{\Gamma}$ through  
    $$\mu'_j:=\max\;\left\{0,\; \mu^*_i+c_i-c_j\right\}, $$
    $\mu'_k:=\mu^*_k$ for all $k\in Q'\setminus\{j\}$, and $\mu'_k:=0$ for all $k\in\bar Q'$.
    Then, we have
    \begin{equation}
    \label{eq:cost_inQ'}
    c_j+\mu'_j\geq c_i+\mu^*_i, \text{ and } c_k+\mu'_k=c_k+\mu^*_k, \forall k\in Q'\setminus\{j\}.
    \end{equation}
    Observe that the values of $c+\mu'$ remain nondecreasing with respect to the index in $Q'$.
    
    Moreover, we obtain
    \begin{align*}
        \bar \Gamma_{Q'}(\mu')
        &=\bar \Gamma_{Q^*}(\mu^*)+\frac{\mu^*_i-\mu'_j}{d}\\
        &=\bar \Gamma_{Q^*}(\mu^*)+\left\{
        \begin{array}{lll}
        &\displaystyle \frac{\mu^*_i}{d}, &\text{ if } c_j>c_i+\mu^*_i,\\
        &\displaystyle\frac{c_j-c_i}{d}, &\text{ otherwise. }
        \end{array}
        \right.   
    \end{align*}
    We now show that any candidate optimal solution $Y'$ to the minimization problem defining $\varphi(Q',\mu')$ has an objective value no smaller than $\varphi(Q^*,\mu^*)$. This can be shown by constructing a feasible solution $Y^*$ to the minimization problem defining $\varphi(Q^*,\mu^*)$ whose objective value is no greater than that of $Y'$. Denote 
    \begin{align*}
        &\Delta:=\sum_{k\in Y'\cap Q'}(c_k+\mu'_k)-\sum_{k\in Y^*\cap Q^*}(c_k+\mu^*_k)\\&\hspace{1.5cm}+\sum_{k\in Y'\cap \bQ'}c_k-\sum_{k\in Y^*\cap \bQ^*}c_k+ d\left(\min\left\{|Y'\cap \bQ'|,\; \bar \Gamma_{Q'}(\mu')\right\}
        -\min\left\{|Y^*\cap \bQ^*|,\; \bar \Gamma_{Q^*}(\mu^*)\right\}\right)
    \end{align*}
    as the difference of their objective value.
    Note that, in what follows, we only need to consider solutions $Y'$ satisfying Remark~\ref{remark:opt_y}. Observe that Remark~\ref{remark:opt_y} implies that $|Y'\cap \{i,j\}|\leq 1$. 
    
    If $i\notin Y'$, then we define $Y^*$ such that $j\notin Y^*$, $i\in Y^*\iff j\in Y'$, and $k\in Y^*\iff k\in Y'$ for each $k\in[n]\setminus\{i,j\}$.  
    It follows from the fact that $\bar \Gamma_{Q'}(\mu')\geq\bar \Gamma_{Q^*}(\mu^*)$ and \eqref{eq:cost_inQ'} that $\Delta\geq 0$. 
    
    If $i\in Y'$ and $j\notin Y'$, 
    we consider the following subcases. When $\bar \Gamma_{Q'}(\mu')\leq |Y'\cap \bQ'|$ and $c_j+\mu'_j=c_i+\mu^*_i$, we use the same construction as above and define $i\notin Y^*$, $j\in Y^*$ and $k\in Y^*\iff k\in Y'$ for each $k\in[n]\setminus\{i,j\}$.
    Then, we have
    \begin{align*}
        \Delta = c_i-c_j+d\cdot(\bar \Gamma_{Q'}(\mu')-\bar \Gamma_{Q^*}(\mu^*))
        =c_i-c_j+c_j-c_i=0.
    \end{align*}
    When $\bar \Gamma_{Q'}(\mu')\leq |Y'\cap \bQ'|$ and $c_j>c_i+\mu^*_i$, 
    since $d\cdot(\bar \Gamma_{Q'}(\mu')-\bar\Gamma_{Q^*}(\mu^*))$ is not necessarily equal to $c_j-c_i$, we consider an alternative construction for $Y^*$ 
    such that $|Y^*\cap Q^*|=|Y'\cap Q'|+1$ and Remark~\ref{remark:opt_y} is satisfied. Let $r$ be the smallest index such that $r\notin Y'$ among the elements of $Q'$ and $s$ be the largest index such that $s\in Y'\cap \bar Q'$. Since $r\leq j$, we have 
     $$
     \sum_{k\in Y'\cap Q'}(c_k+\mu'_k)-\sum_{k\in Y^*\cap Q^*}(c_k+\mu^*_k)\geq -c_i-\mu^*_i.
     $$
     Then, we also have 
     $$
     \sum_{k\in Y'\cap \bar Q'}c_k-\sum_{k\in Y^*\cap\bar Q^*}c_k=c_i \text{ if } s\leq j-1;\; =c_i+c_s-c_j \text{ if } s> j.
     $$
     We thus obtain
    \begin{align*}
        \Delta&\geq -c_i-\mu^*_i+c_i+ d\left(\bar \Gamma_{Q'}(\mu')
        -\min\left\{|Y^*\cap\bar Q^*|,\; \bar\Gamma_{Q^*}(\mu^*)\right\}\right)\\
        &\geq -c_i-\mu^*_i+c_i+ d\left(\bar \Gamma_{Q'}(\mu')-\bar\Gamma_{Q^*}(\mu^*)\right)\\
        &=-c_i-\mu^*_i+c_i+\mu^*_i=0.
    \end{align*}
    When $|Y'\cap\bQ'|<\bar \Gamma_{Q'}(\mu')$, we keep the alternative construction for $Y^*$ 
    such that $|Y^*\cap Q^*|=|Y'\cap Q'|+1$. 
    Then, we have
    \begin{align*}
        \Delta&\geq -c_i-\mu^*_i+c_i+ d\left(|Y'\cap\bQ'|
        -\min\left\{|Y^*\cap\bar Q^*|,\; \bar\Gamma_{Q^*}(\mu^*)\right\}\right)\\
        &\geq -c_i-\mu^*_i+c_i+ d\left(|Y'\cap\bQ'|
        -|Y^*\cap\bar Q^*|\right)\\
        &=-c_i-\mu^*_i+c_i+ d\geq 0.
    \end{align*}
    Consequently,
    $$
    \Phi(Q')=\max_{\mu\in\Xi_{\Gamma}}\varphi(Q',\mu)\geq \varphi(Q',\mu')\geq\varphi(Q^*,\mu^*)=\Phi(Q^*).
    $$    
\end{proof}

\begin{lemma}
\label{lemma:lemma:opt_p}
    There is an optimal solution $Q^*\in\Q$ such that $i<j$ for all $i\in \bQ^*\cap[p]$ and $j\in Q^*\cap[p]$.
\end{lemma}
\begin{proof}
    Denote 
    $$
    Q'=\{i, i+1, \ldots,j-1,q_1,\ldots,q_2\},
    $$
    and
    $$
    \bQ'=\{1,\ldots,i-1,j,\bar{q}_1, \ldots,\bar{q}_2\},
    $$
    where $i,j\in [p]$, $i$ is the first element in $Q'$, and $j$ is the smallest element in $\bQ'$ greater than $i$.
    Define
    $$
    Q^*=Q'\setminus\{i\}\cup \{j\}, \text{ and } \bQ^*=\bQ'\setminus\{j\}\cup \{i\},
    $$    
    namely, 
    $$
    Q^*=\{i+1,\ldots,j-1,j,q_1,\ldots,q_2\},
    $$
    and 
    $$
    \bQ^*=\{1,\ldots,i-1,i,\bar{q}_1, \ldots,\bar{q}_2\}.
    $$
    Note that $j$ and $i$ occupy the same position, namely $i$-th position, in $\bQ'$ and $\bQ^*$, respectively.
    Let $\mu^*\in \argmax_{\mu \in \Xi_{\Gamma}}\varphi(Q^*,\mu)$ satisfying Observation~\ref{obs:opt_Q}, and, without loss of generality, assume that $\mu^*_k=0$ for each $k\in\bar Q^*$. 
    We next show that $\Phi(Q^*)\leq \Phi(Q')$. To do so, we construct a scenario $\mu'\in \Xi_{\Gamma}$ through 
    $$
    \mu'_i:=\min\left\{d,\;\mu^*_j+ c_j-c_i,\;\mu^*_j+d\cdot\bar\Gamma_{Q^*}(\mu^*)\right\},
    $$
    $\mu'_k:=\mu^*_k$ for all $k\in Q'\setminus\{i\}$, and $\mu'_k:=0$ for all $k\in\bar Q'$.
    Observe that $c_i+\mu'_i\leq c_j+\mu^*_j$. However, the values of $c+\mu'$ are not necessarily nondecreasing with respect to the indices in $Q'$, due to the constructed value of $c_i+\mu'_i$.
    Furthermore, $\bar \Gamma_{Q'}(\mu')
        =\bar\Gamma_{Q^*}(\mu^*)+\frac{\mu^*_j-\mu'_i}{d}$. Therefore, whenever $\mu'_i\neq\mu^*_j+d\cdot\bar\Gamma_{Q^*}(\mu^*)$, we have
    \begin{align*}
        \bar \Gamma_{Q'}(\mu')
        &=\bar\Gamma_{Q^*}(\mu^*)-\left\{
        \begin{array}{lll}
        &\displaystyle 1-\frac{\mu^*_j}{d}, &\text{ if } c_j+\mu^*_j>c_i+d,\\
        &\displaystyle\frac{c_j-c_i}{d}, &\text{ otherwise, }
        \end{array}
        \right.   
    \end{align*}
    while $\bar \Gamma_{Q'}(\mu')=0$ if $\mu'_i=\mu^*_j+d\cdot\bar\Gamma_{Q^*}(\mu^*)$. 
    We now show that any candidate optimal solution $Y'$ to the minimization problem defining $\varphi(Q',\mu')$ has an objective value no smaller than $\varphi(Q^*,\mu^*)$.
    This can be shown by constructing a feasible solution $Y^*$ to the minimization problem defining $\varphi(Q^*,\mu^*)$ whose objective value is no greater than that of $Y'$.
    Denote 
    \begin{align*}
        &\Delta:=\sum_{k\in Y'\cap Q'}(c_k+\mu'_k)-\sum_{k\in Y^*\cap Q^*}(c_k+\mu^*_k)\\&\hspace{1.5cm}+\sum_{k\in Y'\cap\bQ'}c_k-\sum_{k\in Y^*\cap\bQ^*}c_k+ d\left(\min\left\{|Y'\cap\bQ'|,\; \bar \Gamma_{Q'}(\mu')\right\}
        -\min\left\{|Y^*\cap\bQ^*|,\; \bar\Gamma_{Q^*}(\mu^*)\right\}\right)
    \end{align*}
    as the difference of their objective value.
    Note that, in what follows, we only need to consider solutions $y'$ satisfying Remark~\ref{remark:opt_y}. 
    Observe that Remark~\ref{remark:opt_y} implies that $|Y'\cap \{i,j\}|\geq 1$. 
        
    We first consider the case $\mu'_i\neq\mu^*_j+d\cdot\bar\Gamma_{Q^*}(\mu^*)$.
    If $i,j \in Y'$, then we define $Y^*:=Y'$, 
    so
    \begin{align*}
        \Delta&= (c_i+\mu'_i-c_j-\mu^*_j)+c_j-c_i+d\left(\min\left\{|Y'\cap\bQ'|,\; \bar \Gamma_{Q'}(\mu')\right\}
        -\min\left\{|Y^*\cap\bQ^*|,\; \bar\Gamma_{Q^*}(\mu^*)\right\}\right)\\
        &\geq\mu'_i-\mu^*_j+d\cdot
        \begin{cases}
        \displaystyle \left(|Y'\cap\bQ'|- |Y^*\cap\bQ^*|\right)= 0, &\text{ if } |Y'\cap\bQ'|\leq \bar\Gamma_{Q'}(\mu'),\\
        \displaystyle \left(\bar\Gamma_{Q'}(\mu')-\bar\Gamma_{Q^*}(\mu^*)\right)=\frac{c_i-c_j}{d}, &\text{ if }|Y'\cap\bQ'|> \bar\Gamma_{Q'}(\mu')\text{ and } c_j+\mu^*_j=c_i+\mu'_i,\\
        \displaystyle \left(\bar\Gamma_{Q'}(\mu')-\bar\Gamma_{Q^*}(\mu^*)\right)= \frac{\mu^*_j}{d}-1, &\text{ if }|Y'\cap\bQ'|> \bar\Gamma_{Q'}(\mu') \text{ and } c_j+\mu^*_j>c_i+d,
        \end{cases}
        \\
        &\geq0.
    \end{align*}
    If $i\notin Y'$ and $j\in Y'$, then we define $Y^*$ such that $i\in Y^*$, $j\notin Y^*$ and $k\in Y^*\iff k\in Y'$ for each $k\in[n]\setminus\{i,j\}$.
    We obtain, by a similar argument to that for the case $i,j\in Y'$, that
    \begin{align*}
        \Delta=c_j-c_i+d\left(\min\left\{|Y'\cap\bQ'|,\; \bar\Gamma_{Q'}(\mu')\right\}
        -\min\left\{|Y^*\cap\bQ^*|,\; \bar\Gamma_{Q^*}(\mu^*)\right\}\right)\geq0.
    \end{align*}
    If $i\in Y'$ and $j\notin Y'$, we consider the following cases. When $c_i+\mu'_i=c_j+\mu^*_j$ and $|Y'\cap\bQ'|\leq\bar\Gamma_{Q'}(\mu')$, we define $Y^*$ such that $j\in Y^*$, $i\notin Y^*$, and $k\in Y^*\iff k\in Y'$ for each $k\in[n]\setminus\{i,j\}$.
    Since we have 
    $$
    \sum_{k\in Y'\cap Q'}(c_k+\mu'_k)=\sum_{k\in Y^*\cap Q^*}(c_k+\mu^*_k),
    $$
    and $k\in Y^*\iff k\in Y'$ for all $k\in \bQ'\cap \bQ^*$, we obtain
    \begin{align*}
        \Delta&=d\left(\min\left\{|Y'\cap\bQ'|,\; \bar\Gamma_{Q'}(\mu')\right\}
        -\min\left\{|Y^*\cap\bQ^*|,\; \bar\Gamma_{Q^*}(\mu^*)\right\}\right)\\
        &= d\left(|Y'\cap\bQ'|-|Y^*\cap\bQ^*|\right)=0.
    \end{align*}
    When $c_i+\mu'_i=c_j+\mu^*_j$ and $|Y'\cap\bQ'|>\bar\Gamma_{Q'}(\mu')$, we define $Y^*$ such that $|Y^*\cap Q^*|=|Y'\cap Q'|-1$ and Remark~\ref{remark:opt_y} is satisfied. Hence, there is $s\in (Q'\cap Q^*)$ with $s>j$ such that $s\in Y'$ and $s\notin Y^*$, and $k\in Y^*\iff k\in Y'$ for all $k\in (Q'\cap Q^*)\setminus\{s\}$. Moreover, let $r\leq i$ be the largest index such that $r\in Y^*\cap \bQ^*$. Then $k\in Y^*\iff k\in Y'$ for all $k\in (\bar Q'\cap \bQ^*)\setminus\{r\}$, and $r\notin Y'$ if $r\neq i$. 
    Therefore,
    \begin{align*}
        \Delta&=c_s+\mu'_s-c_r+    d\left(\bar\Gamma_{Q'}(\mu')
        -\min\left\{|Y^*\cap\bQ^*|,\; \bar\Gamma_{Q^*}(\mu^*)\right\}\right)\\
        &\geq c_s+\mu'_s-c_r+c_i-c_j\geq 0,
    \end{align*}
    where we used $d\cdot(\bar\Gamma_{Q'}(\mu')-\bar\Gamma_{Q^*}(\mu^*))= c_i-c_j$ in the first inequality, and
    $s>j$ and $r\leq i$ in the second one.     
    When $c_j+\mu^*_j>c_i+d$ ($\implies \mu'_i=d$), we also define $Y^*$ such that $|Y^*\cap Q^*|=|Y'\cap Q'|-1$  Defining $r,s$ as above, we have 
    \begin{align*}
        \Delta &=(c_i+d+c_s+\mu'_s-c_j-\mu^*_j)-c_r+    d\left(\min\left\{|Y'\cap\bQ'|,\; \bar\Gamma_{Q'}(\mu')\right\}
        -\min\left\{|Y^*\cap\bQ^*|,\; \bar\Gamma_{Q^*}(\mu^*)\right\}\right)\label{delta}\\
        &\geq\left\{
        \begin{array}{lll}
        &\displaystyle c_i+d-c_r +c_s+\mu'_s-c_j-\mu^*_j-d\geq 0, &\text{ if } |Y'\cap\bQ'|<\bar\Gamma_{Q'}(\mu'),
        \\
        &\displaystyle  c_i+d-c_r +c_s+\mu'_s-c_j-\mu^*_j-(d
        -\mu^*_j)\geq 0, &\text{ if } |Y'\cap\bQ'|\geq\bar\Gamma_{Q'}(\mu'),
        \\
        \end{array}
        \right.   
    \end{align*}
    because $s>j$ and $r\leq i$. 
    
    We now consider the case $\mu'_i=\mu^*_j+d\cdot\bar\Gamma_{Q^*}(\mu^*)$ with $\bar\Gamma_{Q'}(\mu')=0$. 
    We first have $c_i+\mu'_i\leq c_j+\mu^*_j$, 
    which implies that $c_i+d\cdot\bar\Gamma_{Q^*}(\mu^*)\leq c_j$. If $i,j\in Y'$, we define $Y^*:=Y'$. 
    Then, we have 
    \begin{align*}
        \Delta&\geq c_i+\mu^*_j+d\cdot\bar\Gamma_{Q^*}(\mu^*)-c_j-\mu^*_j+c_j-c_i-d\cdot\bar\Gamma_{Q^*}(\mu^*)=0.
    \end{align*}
    If $i\notin Y'$ and $j\in Y'$, then we define $Y^*$ such that $j\notin Y^*$, $i\in Y^*$, and $k\in Y^*\iff k\in Y'$ for each $k\in[n]\setminus\{i,j\}$.
    Then, we have
    \begin{align*}
        \Delta&\geq c_j-c_i-d\cdot\bar\Gamma_{Q^*}(\mu^*)\geq 0.
    \end{align*}
    If $i\in Y'$ and $j\notin Y'$, then we define $Y^*$ such that $|Y^*\cap Q^*|=|Y'\cap Q'|-1$. Then, we have
    \begin{align*}
        \Delta\geq c_i+\mu^*_j+d\cdot\bar\Gamma_{Q^*}(\mu^*)+c_s+\mu'_s-c_j-\mu^*_j
        -c_r-d\cdot\bar\Gamma_{Q^*}(\mu^*)\geq 0,
    \end{align*}
    where $s>j$ and $r\leq i$.

    Consequently,
    $$
    \Phi(Q')=\max_{\mu\in\Xi_{\Gamma}}\varphi(Q',\mu)\geq \varphi(Q',\mu')\geq\varphi(Q^*,\mu^*)=\Phi(Q^*).
    $$    
\end{proof}

\begin{corollary}
\label{cor:optimal_Q}
    There exists an optimal solution $Q^*=\{k,k+1,\ldots,k+q-1\}$ for some $k\in [n-q+1]$.
\end{corollary}
\begin{proof}
    By lemmas~\ref{lemma:opt_bar_p} and~\ref{lemma:lemma:opt_p}.
\end{proof}

    Corollary~\ref{cor:optimal_Q} implies that there are at most $n$ candidate solutions $Q$ whose optimality needs to be verified. While it is known from~\cite{OmerPossRougier2024} that $\Phi(Q)$ can be computed in polynomial time for easy-to-optimize sets $\Y$, for the sake of completeness, we detail this reformulation for the special of interest in what follows. Throughout the following reformulations, we use binary decision variables $y$ rather than their corresponding set representation $Y$, in order to express the problem as a linear program. Therefore, we define $\Y_{IP}:=\set{y\in\{0,1\}^n}{\exists Y\in\mathcal{Y}\text{ s.t. }y_i=1\iff i\in Y}=\set{y\in\{0,1\}^n}{\sum_{i\in[n]}y_i=p}$. First, we derive its  epigraphical reformulation from~\eqref{eq:re_ddid_equal_dev}:
    \begin{align}
        \Phi(Q)=&\max_{\mu,\omega}\; \omega\\
        &\text{ s.t.} \; \omega \leq \min_{y\in \Y_{IP}}\; \sum_{i\in Q}(c_i+\mu_i)y_i+\sum_{i\in \bQ} c_iy_i+ d\cdot \min\left\{\sum_{i\in \bQ}y_i, \; \Gamma-\sum_{i\in Q}\frac{\mu_i}{d}\right\} \label{const:omega_d}\\
        &\quad\quad \mu\in \Xi_{\Gamma}. \label{const:mu_d}
    \end{align}
    Observe that Constraint~\eqref{const:omega_d} can be replaced by the following two constraints:
    \begin{align}
        &\omega \leq \min_{y\in \Y_{IP}}\; \sum_{i\in Q}(c_i+\mu_i)y_i+\sum_{i\in \bQ} c_iy_i+ d\cdot\sum_{i\in \bQ}y_i\label{const:omega1_d}\\
        &\omega \leq \min_{y\in \Y_{IP}}\; \sum_{i\in Q}(c_i+\mu_i)y_i+\sum_{i\in \bQ} c_iy_i+ d\cdot \left(\Gamma-\sum_{i\in Q}\frac{\mu_i}{d}\right)\label{const:omega2_d}
    \end{align}
    Since the matrix defining $\Y_{IP}$ is totally unimodular, the integrality restriction on $y$ can be relaxed in the minimization problems in constraints~\eqref{const:omega1_d} and~\eqref{const:omega2_d}. Dualizing these two linear programs, and removing the maximizations, we obtain
    \begin{align*}
        \Phi(Q)=&\max_{\mu,\omega, \alpha,\alpha',\beta,\beta'}\; \omega\\
        &\text{ s.t.} \; \omega\leq p\cdot \alpha-\sum_{i\in[n]}\beta_i\\
        &\quad\quad \alpha-\beta_i-\mu_i \leq c_i, &\forall i\in Q\\
        &\quad\quad \alpha-\beta_i\leq c_i+d, &\forall i\in \bQ\\
        &\quad\quad \omega\leq p\cdot \alpha'-\sum_{i\in[n]}\beta'_i+d\cdot\Gamma -\sum_{i\in Q}\mu_i\\
        &\quad\quad \alpha'-\beta'_i-\mu_i \leq c_i, &\forall i\in Q\\
        &\quad\quad \alpha'-\beta'_i\leq c_i, &\forall i\in \bQ\\
        &\quad\quad \mu\in\Xi_{\Gamma}\\
        &\quad\quad \beta,\beta'\geq 0.
    \end{align*}
Thus, for any given $Q\in \Q$, the value of $\Phi(Q)$ can be computed by solving linear program defined above, in which there are $O(n)$ variables and $O(n)$ constraints.

\begin{corollary}
\label{cor:complexity_d}
    \textup{\ddid{$\Q^{\SEL}$}{$\Y^{\SEL}$}-\Equal-\textup{\textsf{OPT}}} can be solved by solving at most $n$ linear programs, each with $O(n)$ variables and $O(n)$ constraints.
\end{corollary}

\bibliographystyle{plainurl}
\bibliography{references}

\end{document}